\documentclass{amsart}
\usepackage{graphicx}
\newtheorem{thm}{Theorem}[section]

\newtheorem{prop}[thm]{Proposition}
\theoremstyle{definition}

\theoremstyle{remark}
\newtheorem{rem}[thm]{Remark}
\numberwithin{equation}{section}

\newcommand{\Real}{\mathbb R}

\newcommand{\semi}{\times\!\!\! \mid}

\newcommand{\cS}{\mathcal{S}}

\newcommand{\bZ}{\mathbb{Z}}
\renewcommand{\H}{\mathcal H}
\newcommand{\cen}{\mathfrak Z}
\begin{document}

\title[Generalized Nielsen realization problem]{On the generalized Nielsen realization problem}%
\author{Jonathan Block and Shmuel Weinberger}%
\address{\newline\noindent Department of Mathematics, University of
Pennsylvania\newline Department of Mathematics, University of
Chicago} \email{\newline\noindent blockj@math.upenn.edu,
shmuel@math.uchicago.edu
}%
\begin{abstract}The main goal of this paper is to give the first examples of
equivariant aspherical Poincare complexes, that are not realized
by group actions on closed aspherical manifolds $M$. These will
also provide new counterexamples to the Nielsen realization
problem about lifting homotopy actions of finite groups to honest
group actions. Our examples show that one cannot guarantee that a
given action of a finitely generated group $\pi$ on Euclidean
space extends to an action of $\Pi$, a group containing $\pi$ as a
subgroup of finite index, even when all the torsion of $\Pi$ lives
in $\pi$.
\end{abstract}

\maketitle

\section{Introduction}

Consider an aspherical manifold $M$. Then $\pi_0(\H(M))$, where
$\H(M)$ is the space of self homotopy equivalences of $M$,  is
isomorphic to the group of outer automorphisms of $\pi_1(M)$,
$\mbox{Out}(\pi_1(M))$. The celebrated Borel conjecture,
\cite{FRR} implies that any $\phi\in \H(M)$ is homotopic to a
homeomorphism.  In general, it asserts that homotopy equivalences
(rel boundary, if any) between homotopy equivalent compact
aspherical manifolds are homotopic (rel boundary) to
homeomorphisms.

The Nielsen
realization problem is stated as follows. Given a finite subgroup
$G$ of $\mbox{Out}(\pi_1(M))$, does there exist a group action of
$G$ on $M$ realizing this outer action on $\pi_1(M)$.

In high dimensions, it is easy to give smooth counterexamples to
this using exotic differential structures on the sphere.  Thus, it
makes most sense to consider this problem in topological settings.
We note that, as far as we know, there is no example of
nonrealization even for infinite G.  (However, see \cite{Mor} for
the differentiable failure of this infinite ``Nielsen problem''
for surfaces.)

A first obstruction to $G$ acting on $M$ realizing a given outer action
comes from the nonexistence of certain group extensions.  More
precisely, if the outer action lifts to an actual action, then there is an extension of
groups
\begin{equation}\label{groupextension}
1\to \pi\to \Pi\to G\to 1
\end{equation}
where $\pi=\pi_1(M)$ and the outer action of $G$ on $\pi_1(M)$
arising from the extension is the given one. This condition can be
nontrivial. Raymond and Scott, \cite{RS}, produced examples where
$\pi$ is the fundamental group of a nilmanifold, and for some
cyclic $G$, there exists no such extension (\ref{groupextension}).
However, if the center $\cen(\pi)=0$ there always exists a unique
such an extension, up to isomorphism, \cite{Brown}, Corollary 6.8,
page 106.

Henceforth we assume $\cen(\pi)$ is trivial. Thus in this case
there always exists an extension group $\Pi$, and one reformulates
the Nielsen realization problem and asks if this is enough to
guarantee the existence of an action of $G$ on $M$.

If $\Pi$ is torsion free there is a good conjectural reason to
expect the answer to be positive:
\begin{prop}
If $\Pi$ is torsion free then it is a Poincare duality group if
and only if $\pi$ is. If $B\pi=M$ is a closed manifold of
dimension at least $5$, and the Borel conjecture  holds for $\Pi$
and $\pi$, then $B\Pi$ is a manifold as well and the normal cover
corresponding to $G$ is $M$; thus $M$ has a free $G$ action.

\begin{rem}
We understand the Borel conjecture to assert that if $B\Gamma$ is
any compact manifold with boundary and
$$\phi:(M,\partial M)\to T^n\times(B\Gamma,\partial B\Gamma)$$
is a homotopy equivalence of pairs that is already a homeomorphism
on the boundary, then $\phi$ is homotopic rel boundary to a
homeomorphism. When $B\Gamma$ is a finite complex, this is
well-known to be equivalent to various vanishing statements of
Whitehead groups and isomorphism statements of $L$-theory assembly
maps. In particular, it does not matter which compact manifold
with boundary model of $B\Gamma$ one chooses.
\end{rem}

\end{prop}

\begin{proof} The first statement is Proposition 10.2, page 224 of
\cite{Brown}. As for the second, first observe that $B\Pi$ is a
finite complex by the vanishing of the Wall finiteness obstruction
that lies in the vanishing group $\tilde{K}_0(\bZ\Pi)$.  Now, the
existence of the manifold structure on $B\Pi$ follows from the
theory of the total surgery obstruction:  the obstruction to the
existence of a homology manifold realizing $B\Pi$ lies in a group
which the Borel conjecture asserts is trivial (for this version,
see \cite{BFMW}). This homology manifold is actually a manifold,
because it's covered by one. \end{proof}

\begin{rem} We shall see that the analogue of this proposition for
non-free actions is not true.\end{rem}

One can view the Nielsen problem as one of extending group actions
as follows:  If $\pi$ is the fundamental group of $M$, then $\pi$
naturally acts freely on $\tilde{M}$;  Assuming the extension
$\Pi$ exists, the Nielsen problem asks whether the original $\pi$
action extends to a $\Pi$ action\footnote{Unfortunately, standard
mathematical terminology forces us to overuse the word
``extension".}. (The $\Pi$ action will be free, if and only if
$\Pi$ is torsion free, as in the proposition just discussed.)
Modifying this somewhat, one can ask these extension questions
wherein we demand more on the $\Pi$ action, e.g. that all fixed
sets are empty or contractible (we call this an {\em aspherical
action}, and such an extension of a group action, an {\em
aspherical extension}), cf. e.g. \cite{Lu}, \cite{May}. On the way
to giving our counterexample to Nielsen, we prove the following
theorem which can be thought of as giving a counterexample to
Nielsen realization of free actions on orbifolds.

\begin{thm}\label{maintheorem1}

There is a group extension
$$
1\to \pi\to \Pi\to G\to 1
$$
satisfying the following properties.
\begin{enumerate}
\item Any torsion element in $\Pi$ is in $\pi$, that is $\Pi$ is
relatively torsion free.

 \item $\pi$ is virtually torsion free.

\item $\pi$ acts properly discontinuously and cocompactly on
Euclidean space such that the fixed sets of all finite subgroups
are Euclidean spaces, so $\pi$ is acting aspherically.

 \item The action of $\pi$ does not extend to
one of $\Pi$. In fact, there is no properly discontinuous action
of $\Pi$ on Euclidean space with only contractible fixed-point
sets.

\item There is a properly discontinuous action of $\Pi$ on a
contractible space such that all of the fixed sets of all finite
subgroups are contractible.
\end{enumerate}
\end{thm}

Point (4) above discusses both the statement about free actions on
nonmanifolds and nonfree actions on manifolds.  We give two
constructions. They in fact give a cyclic group of prime order (of
order two for the first construction), $\bZ/p$, which does not act
aspherically on a suitable aspherical manifold.

We also derive \begin{thm}\label{maintheorem2} There is a counter
example to the Nielsen realization problem with group $\bZ/2$ and
centerless fundamental group.
\end{thm}

For a finitely generated discrete group $\Pi$ one can define the
asymptotic homology $HX_*(\Pi)$ of $\Pi$ considered as a metric
space. One has the following dichotomy.
\begin{prop}(\cite{BW}) If $\Pi$ is a group of virtual finite type,
then either $HX_*(\Pi)=\bZ$ for $*=n$ and zero otherwise (which we will call simple) or
$HX_*(\Pi)$ is infinitely generated in some dimension.
\end{prop}

We warn the reader that there are finitely generated groups of
infinite type whose asymptotic homology vanishes in all
dimensions.  For a discrete group $\Pi$ there is a space
$\underline{E\Pi}$, which is universal for proper actions. which
is  unique up to equivariant homotopy equivalence, \cite{Lu} and
\cite{May}. If there is a model for
$\underline{B\Pi}=\underline{E\Pi}/\Pi$ which is a compact
manifold, then the asymptotic homology is simple. It is natural to
ask if this is also sufficient. Our examples answer this as well.
\begin{thm}\label{maintheorem3}
There is a group $\Pi$ of virtual finite type with $HX_*(\Pi)$
simple and which has no proper cocompact action on a contractible
manifold.
\end{thm}

\section{The construction}

For all the theorems above, the constructions are of the following
sort. We will construct $\Pi$ directly via a $\bZ/p$ action on an
aspherical complex, so that properties (2), (4) and (5) either
hold directly by construction, or by computation of a relevant
obstruction. Since this obstruction will vanish on passing to a
finite cover one also obtains the finite index subgroup $\pi$ as
in (3).

We will give two different constructions of such $\Pi$. While they
differ in some details, they both are of the following form. We
will have two aspherical manifolds with boundary $W_1$ and $W_2$,
both boundaries being tori and so that the fundamental group of
the boundary injects. (Or one manifold with two boundary
components.) These manifolds possess $\bZ/p$ actions, but the key
feature is that, while the action on $\partial W_1$ is affine, the
action on $\partial W_2$ is not topologically equivalent to an
affine one. However the actions on the boundaries are
equivariantly homotopy equivalent. Gluing $W_1$ and $W_2$ together
by a homeomorphism homotopic to the equivariant homotopy
equivalence gives a closed manifold $V$ with a homotopy action of
$\bZ/p$ on it, and gluing them together by the equivariant
homotopy equivalence gives the homotopy equivalent complex $X$
with a genuine $\bZ/p$ action. Since the geometric actions on
$W_1$ and $W_2$ are not conjugate, it would seem unlikely that
there would be a corresponding action on the manifold $V$ =
$W_1\cup_\partial W_2$, and showing that will be one of our tasks.
Our debt to \cite{GP} and \cite{J} for inspiration should be
apparent.

Actions on tori with the properties asserted are counterexamples
to the ``equivariant Borel conjecture''. By now, many of these are
known, \cite{CK},\cite{W1}, \cite{W2}, \cite{Sh}. We shall use two
examples: one based on surgery theory (Cappell's Unils) and
another based on embedding theory. The exotic aspherical manifolds
are built by Gromov's hyperbolization, \cite{DJ}, \cite{DJW}.

\subsection{Surgery theory technique}
Consider $\bZ/2$ acting on the torus
$$T=(S^1)^{4n}\times S^1$$
by complex conjugation on the first $4n$ factors and trivially on
the last. The orbifold fundamental group of $T/\bZ/2$ (i.e. the
group of lifts of the action of $\bZ/2$ on the universal cover is
$$ \Gamma =(\bZ^{4n}\semi\bZ/2)\times \bZ$$
Let $a$ be one of Cappell's Unil elements in $L_2(\bZ/2*\bZ/2)$.
Note $\bZ/2*\bZ/2\cong \bZ\semi \bZ/2$. $\Gamma$ retracts onto
$(\bZ\semi\bZ/2)\times\bZ$ and so this class gives rise to a
non-zero class  $\alpha\in L_2(\Gamma)$. So far we have
$(T,\bZ/2)$ with fixed set $F$ a disjoint union of circles. Let
$K$ be the complement of a tubular neighborhood $\mbox{Nbd}(F)$.
Then $\pi_1(K/\bZ/2)\cong\Gamma$. By Wall realization there is a
structure
$$ w(\alpha)\in \cS(K/\bZ/2 \mbox{ rel } \partial)\cong
\cS^{\bZ/2}(K \mbox{ rel } \partial)$$ Now set
$$T'=\mbox{Nbd}(F)\cup w(\alpha)$$
We have thus obtained a new involution on the torus. Moreover $T$
and $T'$ are built equivariantly normally cobordant, call this
normal cobordism $W$. It is not hard to see that the action is not
topologically conjugate to the original affine action, although it
is equivariantly homotopically equivalent to it. (\cite{CK},
\cite{W1}). This can be detected by an element of the isovariant
(that is stratified) structure set in the sense of \cite{W1}.

Now according to \cite{DJW}, we can relatively equivariantly
hyperbolize this normal cobordism $W$ relative to $T\cup T'$ to
get $W_h$, and furthermore, the fundamental groups of the
boundaries still inject into the hyperbolization. The fixed sets
on the boundaries are circles and so the fixed sets in the
cobordism is a surface (of high genus). Now we glue the boundary
components $T$ and $T'$ as described above to get a manifold $V$
and a complex $X$. $X$ is a $\bZ/2$-isovariant aspherical Poincare
complex and $V$ is a manifold with a $\bZ/2$-homotopy action. Let
$$\Pi=\pi_1^{\mbox{orb}}(X)$$
be the orbifold fundamental group of $X$.

Since elements of Unil die on passage to suitable finite covers,
our element $\alpha$ dies when lifted to some finite cover of $T$.
So over $X$ or $V$, the corresponding cover $\hat{X}$ or $\hat{V}$
has an honest manifold structure with an honest $\bZ/2$-action.
Set
$$\pi=\pi_1^{\mbox{orb}}(\hat{V})$$
Then we get
$$1\to\pi\to\Pi\to G\to 1$$
where $G$ is the group of the finite cover. $\pi$ is centerless
since it is an amalgamated free product where one side of the free
product comes from hyperbolization.

We now verify the properties (1)-(5) of Theorem
\ref{maintheorem1}.

\noindent (1) The conjugacy classes of finite order in $\Pi$
correspond to fixed sets in $X$ and thus occur already in $\pi$.

\noindent (2) $\pi$ is virtually torsion free since $\pi\to \bZ/2$
has  torsion free kernel $\pi_1(\hat{X})$ (and $X$ is an
aspherical finite complex).

\noindent (3) We know that $\tilde{X}$ and $\tilde{V}$ are
contractible. Moreover, so are all of their fixed sets.  One can
then cross $X$ and $V$ with $S^1$ (and change $\pi$ to
$\pi\times\bZ$ and $\Pi$ to $\Pi\times\bZ$). This  ensures that
these universal covers are simply connected at infinity and are
thus homeomorphic to Euclidean space.

\noindent (4) We show that $\Pi$ can not act on $\tilde{V}$, as in
the statement of the theorem, with contractible fixed point sets.
If it did, then $\tilde{V}$ is equivariantly homotopy equivalent
to $\tilde{X}$, since $\tilde{V}$ is a model for
$\underline{E\Pi}$, the classifying space for proper actions and
such are unique up to equivariant homotopy equivalence, \cite{Lu}
and \cite{May}. Thus $V$ and its $\bZ/2$-action is equivariantly
homotopy equivalent to $X$ with its action. Note that whenever a
finite group acts on a manifold with manifold fixed sets, then it
also admits such an action with homeomorphic fixed set which is
locally flatly embedded. For a proof of taming theory which
generalizes verbatim to the equivariant situation, see
\cite{Ferry}. Now we can apply a theorem of Browder, \cite{W2},
which says that under a suitable gap and tameness hypotheses, that
isovariant and equivariant homotopy equivalence are the same. So
we conclude that our tamed $\Pi$-space $V$ would be isovariantly
homotopy equivalent to $X$.

Hence it suffices to show that $X$ is not isovariantly homotopy
equivalent to a $\bZ/2$-manifold. Further it therefore suffices to
show that $Y=(X-(X^{\bZ/2}))/(\bZ/2)$ does not have the proper
homotopy type of a manifold. We thus calculate the proper total
surgery obstruction of $Y$. We have the following diagram:

\begin{equation}\begin{array}{ccc}
W_h &\stackrel{\phi }{\longrightarrow } & W \\
    & \stackrel{\tilde{\phi}}\searrow & \downarrow\psi \\
         &      &   T\times I\end{array}
         \end{equation}
All three maps are degree one normal maps.  By \cite{DJ}, $W_h$ is
normally cobordant to $W$ and hence $\phi$ has zero surgery
obstruction. $\psi$ on the other hand has surgery obstruction the
original element $a\in L_2(\Gamma)$.

Now set
$$W_b=W_h\cup_{T\coprod T'} (-W)$$ glueing the boundaries together as before.  But this time we get a
manifold. The surgery obstruction of $W_b\to X$ is still the
original $a$. This obstruction is an element of $L_2(\Gamma,
\pi_1^\infty(Y))$ where of course $\pi_1^\infty(Y)$ is a groupoid
and not a group since $Y$ is not connected at infinity. This maps
to
$$ L_2(\bZ\times(\bZ/2*\bZ/2),\coprod \bZ\times \bZ/2\mbox{'s})$$
We can analyze this by looking at the exact sequence of a pair
$$\cdots \to L_n(\coprod \bZ\times \bZ/2\mbox{'s})\to L_n(\bZ\times(\bZ/2*\bZ/2))
\to L_n(\bZ\times(\bZ/2*\bZ/2),\coprod \bZ\times
\bZ/2\mbox{'s})\cdots$$

\noindent According to Shaneson for any $G$ (ignoring decorations
which we can do since $\bZ/2*\bZ/2$ has vanishing $K$-theory)
$$L_n(\bZ\times G)\cong L_n(G)\times
L_{n-1}(G)$$
and according to Cappell for any $G$ and $H$
$$\tilde{L}(G*H)=\tilde{L}(G)\times \tilde{L}(H)\times Unil(e;
G,H)$$

\noindent  Hence the original element of $Unil$ survives inclusion
into the relative group. Therefore the surgery obstruction of this
normal map is non-zero.

Of course for any other degree one normal map  the same reasoning
shows that  the difference between its surgery obstruction and the
one above lies in the image of the assembly map for
$$H_*(B(\bZ\times(\bZ/2*\bZ/2)),\coprod
B(\bZ\times \bZ/2); \mathbb{L}(e))\to
L_*(\bZ\times(\bZ/2*\bZ/2),\coprod \bZ\times \bZ/2\mbox{'s})
$$
But now, as noted above, the image of this latter group in
$\mbox{Unil}$ is trivial, so we are done.\qed

\proof (of Theorem \ref{maintheorem2}) We begin with the
aspherical manifold $V$ constructed above. In this case set
$\pi=\pi_1 V$. This is centerless as remarked above. $V$ also has
its $\bZ/2$-homotopy action and therefore acts on $\pi$ and $\Pi$
is the semi-direct product. We now argue that the $\pi$-action
does not extend to $\Pi$. This is simply a matter of showing that
any action of $\Pi$ on $V$ automatically has contractible manifold
fixed sets so that we can appeal to the proof of Theorem
\ref{maintheorem1}.

Now, by Smith theory, the fixed set is a $\bZ/2$-homology manifold
homology equivalent mod 2 to $\Real^2$ (by comparison with the
Poincare model $X$.) By \cite{Bredon}, Theorem 16.32, page 388,
for any $p$, any second countable $\bZ/p$-homology manifold of
dimension less than or equal to two is a topological manifold.
Thus, the fixed set is a $2$-manifold which the classification of
surfaces implies that any mod 2 acyclic surface is $\Real^2$. \qed

\begin{rem} Connolly-Davis,\cite{CD}, completed the computation of
$L_n(\bZ/2*\bZ/2,\omega)$ for all $n$ and all orientation
characters $\omega$. As a result, one can modify the above
construction using orientation reversing involutions on tori with
isolated fixed sets, to produce different examples. Given the
calculations of Connolly and Davis, the proof that these examples
work is even more elementary with regard to the verification of
manifoldness of putative fixed sets: the characterization of the
circle is much more straightforward.
\end{rem}

\subsection{Embedding theory technique}
We now give a construction, based on embedding theory, that
suffices for an alternate proof of Theorem \ref{maintheorem1},
which gives examples for $\bZ/p$ for $p$ odd. These are
insufficient for the Nielsen problem since the fixed sets will be
of higher dimension and so we have no way of seeing that they are
automatically manifolds, as in the proof of Theorem
\ref{maintheorem2}.

Let $W_1=T^n\times S^\circ$ where $S^\circ$ is a punctured surface
and $n=2p-4$. Now $\bZ/p$ acts on $W_1$ by permuting the first $p$
circles of $T^n$ leaving the other factors fixed. Let
$F=W_1^{\bZ/p}$. Then $\mbox{dim}(W_1)=2p-2$ and
$\mbox{dim}(F)=p-1$.

We now build a second manifold $W_2$ with a group action by first
producing a new embedding of the fixed set in the boundary torus
$T^n$ using the following general construction, called a finger
move, \cite{Sh}: Let $M^k\subset N^{2k+1}$ be an embedding of
manifolds. Let $[\gamma]\in \pi_1(N)$ be a class represented by a
path $\gamma$ which intersects $M$ only in its two distinct
endpoints, which are assumed to lie in a little ball. Let $R$ be a
regular neighborhood of $\gamma$, a $2k+1$-disk. Then $R\cap
M=D^k\cup D^k$. Move one of the disks $D^k$ along $\gamma$ to have
rel $\partial $ linking number one with the other disk. Remove one
disk of intersection and glue in the other one. We thus arrive at
a new manifold pair $(\mbox{Fing}(N,M,\gamma), M)$ where
$\mbox{Fing}(N,M,\gamma)$ is homeomorphic to $N$ and $M$ is
embedded differently. We can perform the same construction
relative to any finite collection of disjoint curves
$\gamma_1,\cdots, \gamma_k$.

\begin{figure}
\begin{center}
\includegraphics[width=\textwidth]{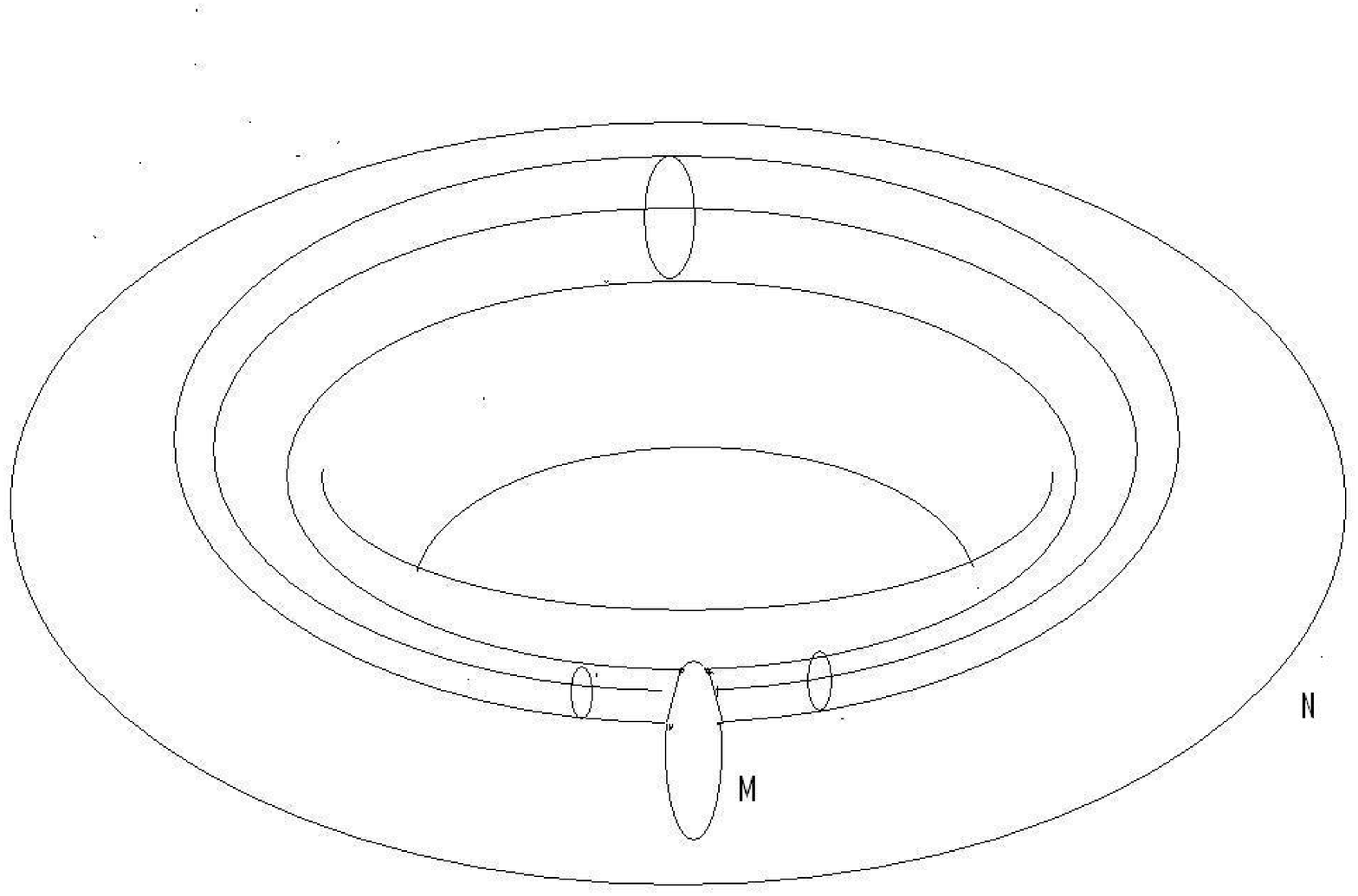}
\caption{A submanifold $M$ of $N$ together with a curve, as data
for a finger move.}
\includegraphics[width=\textwidth]{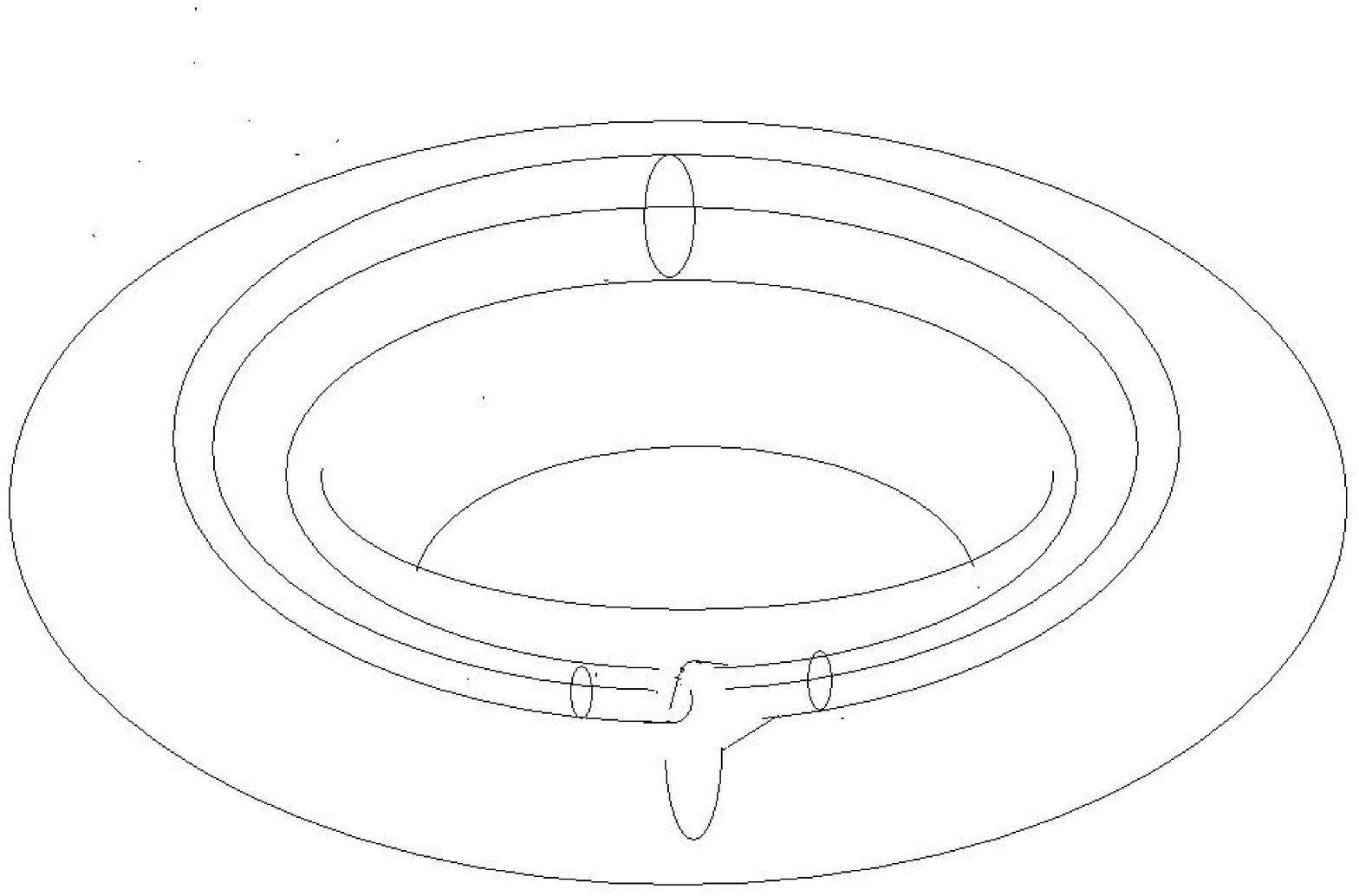}
\caption{The result after the finger move.}
\end{center}
\end{figure}
\newpage
Back to our manifold $W_1$ with its $\bZ/p$-action. Let $\gamma$
be a curve in $\partial W_1=T$ a torus. We may arrange this curve
so that it and all its translates $\gamma, g\gamma,\ldots,
g^{p-1}\gamma$ are disjoint. Now perform the finger move
$W_2=\mbox{Fing}(\partial W_1,\partial W_1\cap F,\gamma, g\gamma,
\ldots, g^{p-1}\gamma)$. We get a new embedding $F'\subset W_2$
and moreover $gF'$ is isotopic to $F'$.

By the main theorem of \cite{W1}, at the cost of repeating all of
these finger moves some number $p^k$ of times, we can find an
equivariantly homotopy equivalent group action on $T$ with fixed
point set $F'$. This action, while a priori only continuous, can
be made PL locally linear (even smooth) and equivariantly
cobordant to the original action on $T$.  This is because
equivariant smoothing theory \cite{LR} and cobordism theory
reduces such problems to the tangent bundle, but \cite{FeW} (see
\cite{FRW} shows that equivariantly homotopy equivalent $G$-tori
have topologically equivalent tangent bundles.

Now we can do our relative hyperbolizations and equivariant
glueing as before to obtain a $\bZ /p$-CW complex unequivariantly
homotopy equivalent to an aspherical manifold W.  We claim that
that this $\bZ/p$-CW complex is not equivariantly homotopy
equivalent to a manifold.  The reason is simple:  the inclusion of
the fixed set $F$ in the $\bZ/p$-CW complex homotopy equivalent to
$W$ is not homotopic to an embedding in $W$.  To check this, we
consider the self intersections of any immersion homotopic to this
inclusion. Note that we are in a non-simply connected situation,
so it is appropriate to use the $\bZ[\pi]$-intersection numbers as
in \cite{Wall}; however, since the subobject $F$ is non-simply
connected, they are not as well defined as in Wall's situation, as
explained in \cite{Sh}. The indeterminacy replaces the $\bZ[\pi]$
by $\bZ[\pi'\backslash \pi/\pi']$ (double cosets) where $\pi'$ is
the fundamental group of $F$, because one can change the path from
basepoint to intersection point either on the way there or on the
way back.

Since we are in the middle dimension, there is a $\bZ$'s worth of
ambiguity, which is reflected in the coefficient of the trivial
double coset $\pi' e \pi' = \pi'$, so we ignore this coefficient.
Of course, the finger move construction gives us a nontrivial
element of $Z[\pi_1( T)/ \pi_1(F \cap T)]$:  this is the usual
relation between linking numbers of chains in a boundary and the
intersection number of bounding cycles.  We only need to see that
nothing is lost on passing to the larger group.  Here we have a
trick available because $\pi_1(F \cap T)$ is normal in $\pi_1(T)$:
the double cosets of $\pi_1(F)$ in $\pi_1(F \cup T) = $(the
group!) $\pi_1(T)/ \pi_1(F \cap  T)$.  Now, general nonsense about
amalgamated free products tells us that $\pi_1(F \cup T)$ injects
into $\pi_1(W)$, so we lose no information at this stage of our
formation of intersection numbers.

Thus, $F$ does not embed in $W$, and therefore neither does any
manifold homotopy equivalent to $F$ in any manifold homotopy
equivalent to W (see e.g. Wall, \cite{Wall}, chapter 11 on
embeddings). A fortoriori, the group action does not exist and our
proof is complete. \qed

\end{document}